\documentclass{LMCS}
\usepackage[latin1]{inputenc}
\usepackage{amssymb,latexsym,amsmath,proof,hyperref}

\def\doi{1 (2:2) 2005}
\lmcsheading%
{\doi}
{7}
{}
{}
{Nov. \hfil \phantom{0}3, 2004}
{Sep. \hfil 26, 2005}
{}

\begin{document}

\title{Internalising modified realisability \\ in constructive type theory}
\author[E.\ Palmgren]{Erik Palmgren}
\address{Department of Mathematics, Uppsala University}
\email{palmgren@math.uu.se}


\begin{abstract} 
A modified realisability interpretation of infinitary logic
is formalised and proved sound in constructive type
theory (CTT). The logic considered subsumes first order logic. The 
interpretation makes it possible to extract programs with simplified
types and to incorporate and reason about them in CTT.
\end{abstract}

\maketitle

\section{Modified realisability}

Modified realisability interpretation is a well-known method for giving 
constructive interpretation of some intuitionistic logical system into 
a  simple type structure \cite{Troelstra:Meta}.
The method is used, for instance, in Minlog and Coq for
extracting programs from proofs (cf.\ \cite{Schwicht:MLCF} and \cite{LetouzeyThesis}). These programs are to a large 
extent free from the computationally irrelevant parts that might be present
in programs arising from direct interpretations into constructive type theory.
The realisability interpretation requires a separate proof
of correctness, which is usually left unformalised.

In this note we present a completely formalised modified realisability 
interpretation carried out in the proof support system 
Agda \cite{CCoquand:Agda}.
We shall here use what is called {\em modified realisability with truth}
which has the property that anything realised is also true in the system 
(Theorem \ref{corr}).
One difference from usual interpretations as in Minlog is that the
logic interpreted goes beyond first order logic: it is a
(constructively) infinitary logic, which arises naturally from 
the type-theoretic notion of universe. Our extension to infinitary
logic seems to be a novel result.

Agda is based on Martin-Löf constructive type theory \cite{ML:AITT} with an infinite
hierarchy of universes $\#0={\sf Set}$, $\#1={\sf Type}$, $\#2={\sf Kind}$,
$\#3$, \ldots. Each of these universes is closed under the formation
of generalised inductive data types. 
We define in Agda an inductive type 
${\sf SP}$ of propositions, so called {\em simple propositions}, by induction:
for each small type $A$ (i.e.\ a member of ${\sf Set}$) an atomic 
proposition ${\sf atom}(A):{\sf SP}$ is introduced; ${\sf SP}$ 
contains $\bot$ and is closed under propositional connectives
($\land$, $\lor$, $\to$)
 and for any small type $A$ and any
propositional function $P:A \to {\sf SP}$ the quantified propositions
$\forall(A,P)$ and $\exists(A,P)$ belong to ${\sf SP}$.
There is an obvious homomorphic embedding ${\sf Tp}$ of ${\sf SP}$ into
the small types defined by ${\sf Tp}(\bot) =\emptyset$, 
${\sf Tp}({\sf atom}(A)) = A$, ${\sf Tp}(P \lor Q) = 
{\sf Tp}(P) + {\sf Tp}(Q)$, ${\sf Tp}(P \land Q) = 
{\sf Tp}(P) \times {\sf Tp}(Q)$, ${\sf Tp}(P \to Q) = 
{\sf Tp}(P) \to {\sf Tp}(Q)$, ${\sf Tp}(\forall(A,P)) = (\Pi x:A){\sf
  Tp}(P(x))$ and ${\sf Tp}(\exists(A,P)) = (\Sigma x:A){\sf Tp}(P(x))$.
We shall sometimes write $(\forall x:A)P(x)$ for $\forall(A,P)$ etc.

The simple propositions may be realised by terms from a
simplified type structure. All atomic propositions will be realised 
by the unique element ${\bf elt}$ of the unit type ${\sf Un}$. 
Define another homomorphism ${\sf Cr}$ (for
{\em crude type}) from  ${\sf SP}$ to
small types by letting
\begin{eqnarray*}
{\sf Cr}(\bot) &=& {\sf Un} \\
{\sf Cr}({\sf atom}(A)) &=& {\sf Un} \\
{\sf Cr}(P \land Q) &=& {\sf Cr}(P)  \times {\sf Cr}(Q) \\
{\sf Cr}(P \lor Q) &=& {\sf Cr}(P)  + {\sf Cr}(Q) \\
{\sf Cr}(P \to Q) &=& {\sf Cr}(P)  \to {\sf Cr}(Q) \\
{\sf Cr}(\forall(A,P)) &=&  (\Pi x:A){\sf Cr}(P(x)) \\
{\sf Cr}(\exists(A,P)) &=&  (\Sigma x:A){\sf Cr}(P(x)). \\
\end{eqnarray*}
The only difference from ${\sf Tp}$ is thus in the translation 
of absurdity and atoms. 
We note that a crude type may still be a dependent type,
if the simple proposition is truly infinitary.
For example, this is the case with ${\sf Cr}(\exists(A,P))$,
if $A=N$ and $P(0) = \top$, $P(S(n)) = Q(n) \land P(n)$.

Another variant of the crude 
type map ${\sf Cr}'$ will be
employed in Theorem \ref{th15} below, which 
is defined as ${\sf Cr}$, except that 
$${\sf Cr}'(\exists(A,P)) =  {\sf Un} + (\Sigma x:A){\sf Cr}'(P(x)).$$
The unit type appearing in the disjoint sum ensures that the type
is never empty, which is crucial for interpreting the full 
absurdity axiom.

The modified realisability ${\sf MR}(S,r)$
of a simple proposition $S:{\sf SP}$ by an element of crude type
$r:{\sf Cr}(S)$ is defined  as a small
proposition (or small type)  by the following
recursion on $S$. (We use the identification of
propositions and types for small types, so that $\land$ and $\lor$
are used interchangeably with $\times$ and $+$, respectively.)

\begin{eqnarray*}
{\sf MR}(\bot,r) &=& \bot \\
{\sf MR}({\sf atom}(P),r) &=& P \\
{\sf MR}(A\land B,r) &=& {\sf MR}(A,r.1) \land {\sf MR}(B,r.2) \\
{\sf MR}(A\lor B,{\sf inl}(s)) &=& {\sf MR}(A,s)  \\
{\sf MR}(A\lor B,{\sf inr}(t)) &=& {\sf MR}(B,t)  \\
{\sf MR}(A\to B,r) &=& ({\sf Tp}(A) \to {\sf Tp}(B))\\ 
& & \land\; 
(\Pi s:{\sf Cr}(A))({\sf MR}(A,s) \to {\sf MR}(B,r(s))) 
  \\
{\sf MR}(\forall(A,P),r) &=& (\Pi x:A){\sf MR}(P(x),r(x)) \\
{\sf MR}(\exists(A,P),r) &=& {\sf MR}(P(r.1),r.2). \\
\end{eqnarray*}
Here $r.1$ and $r.2$ denote the first and second projections.

\begin{rem} The above constructions work in many different 
type-theoretic settings. What is needed is a type universe $U$ closed under 
$\Pi$, $\Sigma$, $+$ and containing basic types ${\sf Un}$ 
and $\emptyset$. Moreover the inductive construction ${\sf SP}_U$
is should be made relative to $U$ instead of ${\sf Set}$. Then 
$${\sf Tp}_U : {\sf SP}_U \to U \qquad {\sf Cr}_U : {\sf SP}_U \to U$$
are defined by recursion on ${\sf SP}_U$ similarly to the above, and so is 
$${\sf MR}_U : (\Pi s :{\sf SP}_U)({\sf Cr}_U(s) \to U).$$
\end{rem}

The following correctness, or conservativity, result states that each simple proposition,
which is realised, is also true in the standard interpretation.

\begin{thm} \label{corr}
For any $S:{\sf SP}$ and $r:{\sf Cr}(S)$, if ${\sf MR}(S,r)$
then ${\sf Tp}(S)$.
\end{thm}
\proof The proof goes by induction on $S$.
For $S=\bot$ or $S={\sf atom}(A)$ the result is immediate. For
$S=A\to B$ we took care to define realisability so that this
is direct as well. Here are two examples of the inductive step. 

Suppose
${\sf MR}(A\lor B,r)$. If $r={\sf inl}(s)$, then ${\sf MR}(A,s)$ is
true. By the inductive hypothesis, we get ${\sf Tp}(A)$
and
hence also ${\sf Tp}(A \lor B)$. The argument for $r={\sf inr}(t)$ is
similar.

Assume ${\sf MR}(\forall(A,P),r)$. Let $a \in A$. Then
${\sf MR}(P(a),r(a))$, and so by the inductive hypothesis ${\sf
  Tp}(P(a))$.
Since $a$ was arbitrary we have actually  ${\sf
  Tp}(\forall(A,P))$. \qed

As a corollary there is an extraction theorem for 
$\forall\exists$-formulae:
\begin{cor} \label{AE} For small types $A$ and $B$ and a simple proposition
$P(x,y)$ where $x:A$ and $y:B$, let $$S =(\forall x:A)
  (\exists y: B)P(x,y).$$
If ${\sf MR}(S, r)$ for some $r$, then there is 
some $f:A \to B$ such that 
${\sf Tp}(P(x,f(x)))$ for all $x:A$.
\end{cor}

Thereby the program $f$ extracted also satisfies its specification
${\sf Tp}(P(x,f(x)))$ within type theory. For $P(x,y) = {\sf
  atom}(R(x,y))$ this is equivalent to $R(x,f(x))$.

\begin{rem}
Note the difference in the $\forall$-case from usual interpretations,
which go 
from theories to theories \cite{Troelstra:Meta}. It is not required that 
${\sf Tp}(\Pi(A,P))$ is added to the condition, since this
follows from the correctness
theorem in the present internalised version.
\end{rem}

We present an intuitionistic infinitary propositional logic
 $IPC^-_{\infty}$ in type theory 
in which quantifiers are
understood as infinitary versions of conjunction and disjunction.
The system has a restriction on the absurdity axiom to atomic
formulae.

$$
\small
\begin{array}{cc}
  A\vdash A \qquad
&\qquad \vcenter{\infer{A\vdash C}{A \vdash B & B \vdash C}} \\ [2ex]
  A \vdash {\sf atom}(P), \rlap{\ for any inhabited $P$}\qquad
&\qquad \\ [2ex]
  A \land B  \vdash A \qquad A \land B  \vdash B\qquad
&\qquad \vcenter{\infer{C\vdash A\land B}{C\vdash A & C\vdash B}} \\ [2ex]
  \bot \vdash {\sf atom}(P) \qquad
&\qquad \\ [2ex]
  A \vdash A \lor B \qquad B \vdash A \lor B \qquad
&\qquad\vcenter{\infer{A \lor B \vdash C}{A \vdash C & B \vdash C}} \\ [2.8ex]
  \vcenter{\infer{A \vdash B\to C}{A \land B \vdash C}} \qquad
&\qquad\vcenter{\infer{A \land B \vdash C}{A \vdash B\to C}} \\ [2.8ex]
  \vcenter{\infer{A\vdash \forall(S,P)}{A \vdash P(t) \quad (t:S)}}\qquad
&\qquad\vcenter{\infer{A\vdash P(t)}{A\vdash \forall(S,P) & t:S}} \\ [2.8ex]
  \vcenter{\infer{\exists(S,P) \vdash A}{P(t) \vdash A \quad (t:S)}}\qquad
&\qquad\vcenter{\infer{P(t) \vdash A}{\exists(S,P) \vdash A \quad t:S}}\\
\end{array}$$

\begin{rem}
Note in particular that the existential quantifier is of the
weak kind, as in first order logic. For $S=\emptyset$ each
$\exists(S,P)$ works as absurdity constant. However, if 
we wish to avoid empty sets as types of realisers, the restricted absurdity
axiom $\bot \vdash {\sf atom(P)}$ should be used. The full
absurdity rule can be derived from the restricted one, 
for those propositions which do not include 
quantification over empty sets. By this procedure we can in principle
extract simply typed programs as in Minlog.
\end{rem}

We say that a sequent $A \vdash B$ is {\em ${\sf MR}$-realised,} 
if there is some $r$ such that ${\sf MR}(A\to B,r)$ is true. A 
rule is {\em realised} if whenever all the sequents above the rule 
bar are realised,  then so is the sequent below the bar.

\begin{thm} \label{th14} The axioms and rules of the system 
$IPC^-_\infty$
are ${\sf MR}$-realised.
\end{thm}

To strengthen the weak absurdity axiom to the full axiom
$$\bot \vdash A$$
where $A:{\sf SP}$ may be arbitrary,
we use the crude type map ${\sf Cr}'$ instead and introduce ${\sf
  MR}'$. This is defined recursively as  ${\sf MR}$ apart from the case for
the existential quantifier:

\begin{eqnarray*}
{\sf MR}'(\exists(S,P),{\sf inl}(s)) &=& \bot\\
{\sf MR}'(\exists(S,P),{\sf inr}(t)) &=& 
 {\sf MR}'(P(t.1),t.2).
\\
\end{eqnarray*}

Theorem \ref{corr} and Corollary \ref{AE} now go through with
${\sf MR}'$ and ${\sf Cr}'$ in place of ${\sf MR}$ and ${\sf Cr}$.

The proof of soundness of the logical rules and axioms is similar as
for Theorem \ref{th14},
with the exception for the verification of the absurdity rule, 
and the left existential rule. This requires a special device. 
Namely a function which to each $P:{\sf SP}$
assigns an element, called ${\sf element}(P)$, of ${\sf Cr}'(P)$ is
necessary.  
This function is defined straightforwardly by recursion on $P$. Some 
key clauses are
\begin{eqnarray*}
{\sf element}(\exists(A,P)) &=& {\sf inl}({\bf elt}) \\
{\sf element}(\forall(A,P)) &=& \lambda x.{\sf element}(P(x)) \\
{\sf element}(A\lor B) &=& {\sf inl}({\sf element}(A)). \\
\end{eqnarray*}
Observe that
no such element need to exist when employing the first 
definition of ${\sf Cr}$, e.g.\ in the case 
${\sf Cr}(\exists(\emptyset,P)) = (\Sigma x:\emptyset){\sf Cr}(P(x))$.

\begin{thm} \label{th15} The axioms and rules of the full system 
$IPC_\infty$ ($IPC^-_\infty$ and the full absurdity axiom)
are ${\sf MR}'$-realised.
\end{thm}


We mention some useful mathematical axioms that are realisable:

\begin{lem} For each propositional function $P:{\mathbb N} \to {\sf
    SP}$ the induction scheme
$$P(0) \land (\forall x:{\mathbb N})[P(x) \to P(S(x))]\to (\forall
  x:{\mathbb N})P(x)$$
 is both ${\sf MR}$-realised and ${\sf MR}'$-realised. 
\end{lem}

\begin{lem} For any binary propositional function  $P:A \times B \to
  {\sf SP}$ the type-theoretic choice principle
$$(\forall x:A)(\exists y:B)P(x,y) \to (\exists g:A\to B)(\forall x:A)P(x,g(x))$$
is ${\sf MR}$-realisable. In case $B$ is inhabited, the principle is
${\sf MR}'$-realisable as well.
\end{lem}
\proof The non-trivial part is to prove the second 
statement. Suppose $b_0:B$ and $r:{\sf Cr}'(S)$ and 
$p: {\sf MR}'(S,r)$, where $S = (\forall x:A)(\exists y:B)\,P(x,y)$.
Define an auxiliary operation $f(x,w):(\Sigma y:B){\sf Cr}'(P(x,y))$
where
$x:A$ and $w: {\sf Cr}'((\exists y:B)P(x,y))$, by cases
\begin{eqnarray*}
f(x, {\sf inl}(u)) &=& \langle b_0, {\sf element}(P(x,b_0)) \rangle \\
f(x, {\sf inr}(y)) &=& y.
\end{eqnarray*}
The realiser $k$ for the implication is now given by
$$k(r) = \langle \lambda x. f(x,r(x)).1, \lambda x. f(x,r(x)).2
\rangle$$
To prove it is a realiser, use $\bot$-elimination for the case 
$r(x) = {\sf inl}(u)$.\qed

The following result is often useful to verify realisability.

\begin{lem} If the ${\sf Tp}$-translation of the proposition
$$(\forall x_1:A_1)\cdots (\forall x_n:A_n)[Q(x_1,\ldots,x_n) \to
P(x_1,\ldots,x_n)] $$
is true and $P$ is atomic or $\bot$, then the proposition
 is ${\sf MR}$-realised as well as ${\sf MR}'$-realised. 
\end{lem}
\proof The realising function is trivial for such
a proposition: $(\lambda x_1)\cdots (\lambda x_n)(\lambda r){\bf
  elt}$. Theorem \ref{corr}, a special property of modified realisability with
truth, is necessary here.\qed

Many stronger ``transfer principles'' are possible to establish. See 
\cite{BBS:Refined} for further results and references.

\section{An Example}

We test the formalisation and extraction procedure on a simple example,
which is due to Berger and Schwichtenberg.
The extracted  function computes Fibonacci numbers efficiently by ``memoization.''

A binary predicate $G$ on natural numbers is given. From
the axioms
\begin{itemize}
\item[(Ax1)] $G(0,0)$
\item[(Ax2)] $G(1,1)$
\item[(Ax3)] $(\forall m, k, \ell)[G(m,k) \land G(S(m),\ell) \to G(S(S(m)),k+\ell)].$
\end{itemize}
one derives by induction and intuitionistic logic the proposition
\begin{itemize}
\item[(P)] $(\forall x)(\exists k,\ell)G(x,k) \land G(S(x),\ell)$.
\end{itemize}
Thus there is some realiser $f$ so that 
$${\sf MR}({\rm Ax1}\, \&\, {\rm Ax2}\, \&\, {\rm Ax3} \vdash
{\rm P},f).$$ The extracted program $p$ 
(which is {\tt fib\_prog} in the Appendix) 
for computing the Fibonacci sequence is then given by
$$p(x) = f(nc,x).1$$
where $nc$ ({\tt nocontent} in the Appendix)
is the trivial realiser for ${\rm Ax1}\, \&\, {\rm Ax2}\, \&\, {\rm Ax3}$.
After a normalisation process one gets the program:
\begin{verbatim}
p x =
 
(case x of {
    (zero) -> t;
    (succ x') -> h x' g (rec
                           (\(z::Nat) -> C)
                           x'
                           t
                           (\(x''::Nat) -> \(y::C) -> h x'' g y));}).1
 
\end{verbatim}
%
%
where
\begin{verbatim}

      C = Sigma Nat (\(k::Nat) -> Sigma Nat (\(l::Nat) -> Unit))

h v p q = <q.2.1;
                 <case q.2.1 of {(zero) -> q.1;
                                 (succ u) -> succ (q.1 + u);
                                }
                       ;<q.2.2.2; e>>>
  
      t = <zero; <succ zero; <e;e>>>

      g = \(x,y,z::Nat) -> \(h,j::Unit) -> e

      e = elt@_
\end{verbatim}

\begin{rem} Note that all truly dependent types have disappeared.
The type $C$ is really the type  ${\sf N}\times ({\sf N} \times {\sf Un})$.

The normalised program has been computed using the partial
normalisation procedure of Agda on selected subexpressions, and was
thus
not completely automatic. We
also introduced the abbreviations $C, h,t,g,e$ by hand.
Some syntactical sugar for lambda expressions and pairs is added.
\end{rem}

\section{The formalisation}

The formalisation have been carried out in 
Agda/IAgda (version 2003-08-09) with the aid of the 
graphical user interface Alfa.
The relevant files are available at the URL
\begin{center}
\texttt{www.math.uu.se/\~{}palmgren/modif}
\end{center}

\nocite{BBSS:PfThyWo}
\bibliographystyle{amsalpha}
\bibliography{ModBiblio}

\end{document}